\documentclass[11pt,reqno]{amsart}
\usepackage{graphicx}
\usepackage{amsfonts}
\usepackage{mathrsfs}
\usepackage{times}
\usepackage{mathptmx}
\usepackage{amsmath}
\usepackage{txfonts}
\usepackage{amssymb}
\usepackage{amsmath,amsthm}
\usepackage{amsfonts}

\allowdisplaybreaks[4] \numberwithin{equation}{section}

\newtheorem{thm}{Theorem}

\newtheorem{cor}{Corollary}
\newtheorem{lem}{Lemma}
\newtheorem{rem}{Remark}
\newtheorem{pro}{Proposition}

\def\squarebox#1{\hbox to #1{\hfill\vbox to #1{\vfill}}}

\newcommand{\qbinomial}[3]{\mbox{$
\biggl[ 
\begin{array}{c}
#1\\
 #2
\end{array}\biggr]_{
\!{#3}}$}}

\begin{document}

\begin{center}
\textbf{\large A Note on Generalized $q$-Difference Equations 
and Their Applications\\[1mm] 
Involving $q$-Hypergeometric Functions}\\[4mm]

\textbf{H. M. Srivastava$^{1,2,3,\ast}$, Jian Cao$^{4}$ 
and Sama Arjika$^{5,6}$}\\[2mm]

$^{1}$Department of Mathematics and Statistics, University of Victoria,\\ 
Victoria, British Columbia $V8W \; 3R4$, Canada\\[1mm]

$^{2}$Department of Medical Research, 
China Medical University Hospital, \\
China Medical University, 
Taichung $40402$, Taiwan, Republic of China\\[1mm]

$^{3}$Department of Mathematics and Informatics, Azerbaijan University,\\  
71 Jeyhun Hajibeyli Street, AZ$1007$ Baku, Azerbaijan\\[1mm]

\textbf{E-Mail: harimsri@math.uvic.ca}\\[1mm]
 
{\bf $^{\ast}$Corresponding Author}\\[2mm]

$^{4}$Department of Mathematics, Hangzhou Normal University,\\ 
Hangzhou City 311121, Zhejiang Province, People's Republic of China\\[1mm]

\textbf{E-Mail: 21caojian@hznu.edu.cn}\\[2mm]

$^{5}$Department of Mathematics and Informatics, University of Agadez, Niger\\[1mm]

$^{6}$International Chair of Mathematical Physics and Applications (ICMPA-UNESCO Chair),\\  
University of Abomey-Calavi, Post Box 072, Cotonou 50, Republic of Benin\\[1mm]

{\bf E-Mail: rjksama2008@gmail.com}
\end{center}

\bigskip

\begin{center}
\textbf{Abstract}
\end{center}
\begin{quotation}
In this paper, we use two $q$-operators $\mathbb{T}(a,b,c,d,e,yD_x)$ and \break 
$\mathbb{E}(a,b,c,d,e,y\theta_x)$ to derive two potentially useful 
generalizations of the $q$-binomial theorem, a set of two extensions of 
the $q$-Chu-Vandermonde summation formula and two new generalizations 
of the Andrews-Askey integral by means of the $q$-difference equations.
We also briefly describe relevant connections of various special cases 
and consequences of our main results with a number of known results. 
\end{quotation}

\medskip


$q$-Series and $q$-integral identities; 
$q$-Difference equations; Sears transformation.

\section{\bf Introduction, Definitions and Preliminaries}

Throughout this paper, we refer to \cite{GasparRahman} for definitions 
and notations. We also suppose that $0<q<1$. For complex numbers $a$, 
the $q$-shifted factorials are defined by
\begin{align}
(a;q)_0:=1,\quad  (a;q)_{n} =\prod_{k=0}^{n-1} (1-aq^k)   
\quad \text{and} \quad (a;q)_{\infty}:=\prod_{k=0}^{\infty}(1-aq^{k}),
\end{align}
where (see, for example, \cite{GasparRahman} and \cite{Slater})
$$(a;q)_n=\frac{(a;q)_\infty}{(aq^n;q)_\infty} \qquad \text{and} \qquad
(a;q)_{n+m}=(a;q)_n(aq^n;q)_m$$
and 
$$\left(\frac{q}{a};q\right)_n=(-a)^{-n}\;q^{({}^{n+1}_{\,\,\,2})}
\frac{(aq^{-n};q)_\infty}{(a;q)_\infty}.$$
We adopt the following notation:  
$$(a_1,a_2, \cdots, a_r;q)_m=(a_1;q)_m (a_2;q)_m\cdots(a_r;q)_m
\qquad (m\in \mathbb{N}:=\{1,2,3,\cdots\}).$$
Also, for $m$ large, we have
$$(a_1,a_2, \cdots, a_r;q)_\infty=(a_1;q)_\infty 
(a_2;q)_\infty\cdots(a_r;q)_\infty.$$

The $q$-binomial coefficient is defined by
\begin{equation}
\begin{bmatrix}
n \\
k \\
\end{bmatrix}_q:=\frac{(q;q)_n}{(q;q)_k(q;q)_{n-k}}.
\end{equation}

The basic (or $q$-) hypergeometric function 
of the variable $z$ and with $\mathfrak{r}$ numerator 
and $\mathfrak{s}$ denominator parameters 
is defined as follows (see, for details, the monographs by 
Slater \cite[Chapter 3]{Slater} 
and by Srivastava and Karlsson 
\cite[p. 347, Eq. (272)]{SrivastavaKarlsson}; 
see also \cite{HMS-IMAJAM1983-1984} and \cite{Koekock}): 

$${}_{\mathfrak r}\Phi_{\mathfrak s}\left[
\begin{array}{rr}
a_1, a_2,\cdots, a_{\mathfrak r};\\
\\
b_1,b_2,\cdots,b_{\mathfrak s};
\end{array}\,
q;z\right]
:=\sum_{n=0}^\infty\Big[(-1)^n \; 
q^{\binom{n}{2}}\Big]^{1+{\mathfrak s}-{\mathfrak r }}
\,\frac{(a_1, a_2,\cdots, a_{\mathfrak r};q)_n}
{(b_1,b_2,\cdots,b_{\mathfrak s};q)_n}
\; \frac{z^n}{(q;q)_n},$$
where $q\neq 0$ when ${\mathfrak r }>{\mathfrak s}+1$. 
We also note that

$${}_{\mathfrak r+1}\Phi_{\mathfrak r}\left[
\begin{array}{rr}
a_1, a_2,\cdots, a_{\mathfrak r+1}\\
\\
b_1,b_2,\cdots,b_{\mathfrak r };
\end{array}\,
q;z\right]
=\sum_{n=0}^\infty \frac{(a_1, a_2,\cdots, a_{\mathfrak r+1};q)_n}
{(b_1,b_2,\cdots,b_{\mathfrak r};q)_n}\;\frac{ z^n}{(q;q)_n}.$$\par
We remark in passing that, in a recently-published 
survey-cum-expository review article, the so-called $(p,q)$-calculus 
was exposed to be a rather trivial and inconsequential variation of 
the classical $q$-calculus, the additional parameter $p$ being redundant 
or superfluous (see, for details, \cite[p. 340]{HMS-ISTT2020}).
 
Basic (or $q$-) 
series and basic (or $q$-) polynomials, especially
the basic (or $q$-) hypergeometric functions and basic 
(or $q$-) hypergeometric polynomials, are 
applicable particularly in several areas of Number Theory 
such as the Theory of Partitions and  
are useful also in a wide variety 
of fields including, for example, Combinatorial Analysis,
Finite Vector Spaces, Lie Theory, Particle Physics, Non-Linear
Electric Circuit Theory, Mechanical Engineering, Theory of
Heat Conduction, Quantum Mechanics, Cosmology, and Statistics
(see also \cite[pp. 350--351]{SrivastavaKarlsson} and the references cited
therein). Here, in our present investigation, we are mainly concerned  
with the Cauchy polynomials $p_n(x,y)$ as given
below (see \cite{Chen2003} and \cite{GasparRahman}):
\begin{equation}
\label{def}
p_n(x,y):=(x-y)(x- qy)\cdots (x-q^{n-1}y)
=\left(\frac{y}{x};q\right)_n\,x^n, 
\end{equation}
together with the following Srivastava-Agarwal 
type generating function 
(see also \cite{Cao-Srivastava2013}):  
\begin{equation}
\label{Srivas}
\sum_{n=0}^\infty  
p_n (x,y)\;\frac{(\lambda;q)_n\,t^n}{(q;q)_n}
={}_{2}\Phi_1\left[
\begin{array}{rr}
\lambda,\frac{y}{x};\\
\\
0;
\end{array}\,
q;  xt\right].
\end{equation}

In particular, for $\lambda=0$ in (\ref{Srivas}), we get the 
following simpler generating function \cite{Chen2003}:
\begin{equation}
\label{gener}
\sum_{n=0}^{\infty} p_n(x,y)\;
\frac{t^n }{(q;q)_n} = 
\frac{(yt;q)_\infty}{(xt;q)_\infty}.
\end{equation}

The generating function (\ref{gener}) 
is also the homogeneous version
of the Cauchy identity or the following 
$q$-binomial theorem (see, for example, \cite{GasparRahman}, 
\cite{Slater} and \cite{SrivastavaKarlsson}):
\begin{equation}
\label{putt}
\sum_{k=0}^{\infty} 
\frac{(a;q)_k }{(q;q)_k}\;z^{k}={}_{1}\Phi_0\left[
\begin{array}{rr}
a;\\
\\
\overline{\hspace{3mm}}\,;
\end{array} \,
q;z\right]=\frac{(az;q)_\infty}{(z;q)_\infty}\qquad (|z|<1). 
\end{equation}
Upon further setting $a=0$, this last relation (\ref{putt}) 
becomes Euler's identity  
(see, for example, \cite{GasparRahman}):
\begin{equation}
\label{q-expo-alpha}
\sum_{k=0}^{\infty} \frac{z^{k}}{(q;q)_k}=\frac{1}{(z;q)_\infty}
\qquad (|z|<1)
\end{equation}
and its inverse relation given below \cite{GasparRahman}:
\begin{equation}
\label{q-Expo-alpha}
 \sum_{k=0}^{\infty}\frac{(-1)^k
}{(q;q)_k}\; q^{\binom{k}{2}}\,z^{k}=(z;q)_\infty.
\end{equation}
Based upon the $q$-binomial theorem (\ref{putt}) and Heine's
transformations, Srivastava {\it et al.} \cite{HMS-C-W2020} established 
a set of two presumably new theta-function identities 
(see, for details, \cite{HMS-C-W2020}).

The following usual $q$-difference operators are defined by 
\cite{Liu97,SrivastaAbdlhusein,Saadsukhi}
\begin{equation} 
\label{deffd}
{D}_a\big\{f(a)\big\}:=\frac{f(a)-f( q a)}{ a}, \quad  
{\theta}_{a}=  \big\{f(a)\}:=\frac{f(q^{-1}a)-f( a)}{q^{-1}a},
\end{equation}
and their Leibniz rules are given by (see \cite{Roman1982})
\begin{align}
\label{Lieb}
D_a^n\left\{f(a)g(a)\right\}=\sum_{k=0}^n 
\qbinomial{n}{k}{q} q^{k(k-n)} D_a^k\left\{f(a)\right\} 
D_a^{n-k}\left\{g(q^ka)\right\}
\end{align}
and
\begin{align}
\theta_a^n\left\{ f(a)g(a)\right\}=\sum_{k=0}^n \qbinomial{n}{k}{q}   
\theta_a^k\left\{f(a)\right\} \theta_a^{n-k}\left\{g(q^{-k}a)\right\},
\label{Lieb2}
\end{align}
respectively. Here, and in what follows, 
$D_a^0$  and $\theta_a^0$ are understood as the identity operators.

Recently, Chen and Liu \cite{Liu97,Liu98} constructed the following 
pair of augmentation
operators, which is of great significance for deriving identities 
by applying its various special cases:
\begin{equation}
\mathbb{T}(bD_x)=\sum_{n=0}^\infty\frac{(bD_x)^n}{(q;q)_n}\qquad\text{and}
\qquad  
\mathbb{E}(bD_x)=\sum_{n=0}^\infty\frac{(b\theta_x)^n}{(q;q)_n}.
\end{equation}
Subsequently, Chen and Gu \cite{CHEN2008} defined 
the Cauchy augmentation operators as follows: 
\begin{equation}
\mathbb{T}(a, bD_x)=\sum_{n=0}^\infty\frac{(a;q)_n}{(q;q)_n}\,
(bD_x)^n
\end{equation}
and
\begin{equation} 
\mathbb{E}(a,bD_x)=\sum_{n=0}^\infty\;\frac{(b;q)_n}{(q;q)_n}\;
(-b\theta_x)^n.
\end{equation}
On the other hand, Fang \cite{Fang2010} 
and Zhang and Wang \cite{Zhang2010} 
considered the following finite generalized $q$-exponential 
operators with two parameters:
\begin{equation}
\mathbb{T}\left[\begin{array}{c}q^{-N},w\\
v\end{array}\Big|q;tD_x\right] =\sum_{n=0}^N\;\frac{(q^{-N},w;q)_n}
{(v,q;q)_n}\,(tD_x)^n
\end{equation}
and
\begin{equation}
\mathbb{E}\left[\begin{array}{c}q^{-N},w\\
v\end{array}\Big|q;t\theta_x\right] 
=\sum_{n=0}^N\;\frac{(q^{-N},w;q)_n}{(v,q;q)_n}\,(t\theta_x)^n.
\end{equation}
Moreover, Li and Tan \cite{LiTan2016} constructed two generalized 
$q$-exponential operators with three parameters as follows:

\begin{equation}
\mathbb{T}\left[\begin{array}{c}u,v\\
w\end{array}\Big|q;tD_x\right] =\sum_{n=0}^\infty\;
\frac{(u,v;q)_n}{(w,q;q)_n}
\,(tD_x)^n
\end{equation}
and
\begin{equation}
\mathbb{E}\left[\begin{array}{c}u,v\\
w\end{array}\Big|q;t\theta_x\right] 
=\sum_{n=0}^\infty\;\frac{(u,v;q)_n}{(w,q;q)_n}\,
(t\theta_x)^n.
\end{equation}

Finally, we recall that Cao {\it et al.} \cite{JianCaoArjika} constructed   
the following $q$-operators:  
\begin{align}
\mathbb{T}(a,b,c,d,e,yD_x)=\sum_{n=0}^\infty\;\frac{(a,b,c;q)_n}{(q,d,e;q)_n}
\;(yD_x)^n
\label{1.9}
\end{align}
and
\begin{align}
\mathbb{E}(a,b,c,d,e,y\theta_x)=\sum_{n=0}^\infty\;\frac{(-1)^nq^{n\choose2}
(a,b,c;q)_n}{(q,d,e;q)_n}\;(y\theta_x)^n \label{1.10}
\end{align}
and thereby generalized Arjika's results in \cite{Arjika2020} by 
using the $q$-difference equations (see, for details, \cite{JianCaoArjika}). 

We remark that the $q$-operator (\ref{1.9}) is a particular case of  
the homogeneous $q$-difference operator $\mathbb{T}({\bf a},{\bf b},cD_{x})$  
(see \cite{HMS-Sama2020}) by taking  
$${\bf a} =(a,b,c),\quad {\bf b} =(d,e) \qquad \text{and} \qquad c=y.$$  
Furthermore, for $b=c=d=e=0$, the $q$-operator (\ref{1.10}) reduces to the operator
$\widetilde{L}(a,y;\theta_{x})$ which was investigated by Srivastava 
{\it et al.} \cite{6}.

\begin{pro}{\rm (see \cite[Theorems 3]{JianCaoArjika})}\label{thm2}
Let $f(a,b,c,d,e,x,y)$ be a seven-variable analytic function in a 
neighborhood of $(a,b,c,d,e,x,y)=(0,0,0,0,0,0,0)\in\mathbb{C}^7$.\\

\noindent
{\rm (I)} If $f(a,b,c,d,e,x,y)$ satisfies 
the following difference equation$:$
\begin{align}\label{thm2_1}
&x\big\{f(a,b,c,d,e,x,y)-f(a,b,c,d,e,x,yq)\notag \\
&\qquad \qquad\quad \quad-(d+e)q^{-1}\;
[f(a,b,c,d,e,x,yq)-f(a,b,c,d,e,x,yq^2)]
\nonumber\\
&\qquad \qquad\quad \quad +deq^{-2}\;[f(a,b,c,d,e,x,yq^2)-f(a,b,c,d,e,x,yq^3)]\big\}
\nonumber\\
&\qquad \quad =y\big\{[f(a,b,c,d,e,x,y)-f(a,b,c,d,e,xq,y)]\notag \\
&\qquad \qquad\quad \quad-(a+b+c)[f(a,b,c,d,e,x,yq)-f(a,b,c,d,e,xq,yq)]
\nonumber\\
&\qquad\qquad\quad\quad +(ab+ac+bc)[f(a,b,c,d,e,x,yq^2)-f(a,b,c,d,e,xq,yq^2)]
\nonumber\\
&\qquad\qquad\quad\quad-abc[f(a,b,c,d,e,x,yq^3)-f(a,b,c,d,e,xq,yq^3)]\big\},
\end{align}
then 
\begin{align}\label{thm2_1.2}
f(a,b,c,d,e,x,y)=\mathbb{T}(a,b,c,d,e,yD_x)\{f(a,b,c,d,e,x,0)\}.
\end{align}

\noindent
{\rm (II)} If $f(a,b,c,d,e,x,y)$ satisfies the following difference equation$:$
\begin{align}\label{thm2_2}
&x\big\{f(a,b,c,d,e,xq,y)-f(a,b,c,d,e,xq,yq)\notag \\
&\qquad \qquad\quad \quad-(d+e)q^{-1}\;
[f(a,b,c,d,e,xq,yq)-f(a,b,c,d,e,xq,yq^2)]\nonumber\\
&\quad\quad\quad\quad +deq^{-2}\;
[f(a,b,c,d,e,xq,yq^2)-f(a,b,c,d,e,xq,yq^3)]\big\}
\nonumber\\
&\qquad\quad =y\big\{[f(a,b,c,d,e,xq,yq)-f(a,b,c,d,e,x,yq)]\notag \\
&\qquad \qquad\quad \quad-(a+b+c)
[f(a,b,c,d,e,xq,yq^2)-f(a,b,c,d,e,x,yq^2)]\nonumber\\
&\qquad \qquad\quad\quad +(ab+ac+bc)[f(a,b,c,d,e,xq,yq^3)-f(a,b,c,d,e,x,yq^3)]
\nonumber\\
&\qquad\qquad\quad\quad -abc[f(a,b,c,d,e,xq,yq^4)-f(a,b,c,d,e,x,yq^4)]\big\},
\end{align}
then 
\begin{align}
f(a,b,c,d,e,x,y)=\mathbb{E}(a,b,c,d,e,y\theta_x)\{f(a,b,c,d,e,x,0)\}.
\end{align}
\end{pro}

Liu \cite{Liu2010,Liu2011} initiated the method based upon 
$q$-difference equations and deduced several results involving 
Bailey's ${}_6\psi_6,\,q$-Mehler formulas 
for the Rogers-Szeg\"{o} polynomials 
and $q$-integral version of the Sears transformation.

\begin{lem}
\label{MALM}
Each of the following $q$-identities holds true$:$
\begin{align}\label{id1}
{D}_a^k \left\{\frac{1}{(as;q)_\infty}\right\}
=\frac{s^k}{(as;q)_\infty},
\end{align}
\begin{align}\label{id2}
\theta_a^k \left\{\frac{1}{(as;q)_\infty}\right\}
=\frac{s^kq^{-\binom{k}{2}}}{(asq^{-k};q)_\infty},
\end{align}
\begin{align}\label{id3}
{D}_a^k \left\{{(as;q)_\infty}\right\}
=(-s)^k \;q^{\binom{k}{2}}\; 
(asq^k;q)_\infty,
\end{align}
\begin{align}\label{id4}
\theta_a^k \left\{{(as;q)_\infty}\right\}
=(- s)^k \;(as;q)_\infty,\\
D_a^n\left\{\frac{(as;q)_\infty}{(a\omega;q)_\infty}\right\}
=\omega^n\; \frac{\left(\frac{s}{\omega};q\right)_n}
{(as;q)_n}\;\frac{(as;q)_\infty}{(a\omega;q)_\infty}
\end{align}
and

\begin{align}\label{aberll} 
\theta_a^n\left\{\frac{(as;q)_\infty}{(a\omega;q)_\infty}\right\}  
=\left(-\frac{q}{a}\right)^n\; \frac{\left(\frac{s}{\omega};q\right)_n}
{\left(\frac{q}{a\omega};q\right)_n}   
\; \frac{(as;q)_\infty}{(a\omega;q)_\infty}. 
\end{align} 
\end{lem}

We now state and prove the $q$-difference formulas as
Theorem \ref{dsdpqro1} below.
 
\begin{thm} \label{dsdpqro1} 
Each of the following assertions holds true$:$  
\begin{align}\label{gLEM}
&\mathbb{T}(r,f,g,v,w,uD_a)
\left\{\frac{(as;q)_\infty}{(az,at;q)_\infty}\right\} \notag \\
&\qquad \quad =\frac{(as;q)_\infty}{(az,at;q)_\infty}\;
\sum_{k=0}^\infty\;
\frac{\left(r,f,g,\frac{s}{z},at;q\right)_{k}(zu)^k}
{(v,w,as,q;q)_{k}}\;
{}_3\Phi_2\left[\begin{matrix}
\begin{array}{rrr}
rq^k,fq^k,gq^k;\\\\
vq^k,wq^k;
\end{array}
\end{matrix} q; ut\right]
\end{align}
and
 
\begin{align}\label{gdLEM} 
&\mathbb{E}(r,f,g,v,w,u\theta_a)\left\{\frac{(az,at;q)_\infty}
{(as;q)_\infty}\right\} \notag \\
&\qquad =\frac{(az,at;q)_\infty}{(as;q)_\infty}\;
\sum_{k=0}^\infty 
\; \frac{\left(r,f,g,\frac{z}{s},\frac{q}{at};q\right)_k\,
(-ut)^k}{\left(v,w,\frac{q}{as},q;q\right)_k}\;
{}_3\Phi_3\left[\begin{matrix}
\begin{array}{rrr}
rq^k,fq^k,gq^k;\\\\
vq^k,wq^k,0;
\end{array}
\end{matrix} q; -ut\right],
\end{align}
provided that $\max\left\{|az|,|as|,|at|,|ut|\right\}<1.$
\end{thm}

\begin{proof} 
By the means of the definitions (\ref{1.9}) and (\ref{1.10}) of the
operators $\mathbb{T}(r,f,g,v,w,uD_a)$ 
and $\mathbb{E}(r,f,g,v,w,u\theta_a)$
and the Leibniz rules (\ref{Lieb}) and (\ref{Lieb2}), we observe that 
\begin{align}\label{proof3}
&\mathbb{T}(r,f,g,v,w,uD_a)\left\{\frac{(as;q)_\infty}
{(az,at;q)_\infty}\right\}
\notag \\
&\qquad =\sum_{n=0}^\infty \frac{(r,f,g;q)_nu^n}{(v,w,q;q)_n} 
D_a^n\left\{\frac{(as;q)_\infty}
{(az,at;q)_\infty}\right\}\nonumber\\
&\qquad=\sum_{n=0}^\infty \frac{(r,f,g;q)_nu^n}{(v,w,q;q)_n}\; 
\sum_{k=0}^n \qbinomial{n}{k}{q}\; q^{k(k-n)}\; 
D_a^k\left\{\frac{(as;q)_\infty}
{(az;q)_\infty}\right\} D_a^{n-k}
\left\{\frac{1}{(atq^{k};q)_\infty}\right\}\notag \\
&\qquad=\sum_{n=0}^\infty \frac{(r,f,g;q)_nu^n}{(v,w,q;q)_n} \;
\sum_{k=0}^n \qbinomial{n}{k}{q}\;q^{k(k-n)}
\; \frac{\left(\frac{s}{z};q\right)_k\,z^k}{(as;q)_k}\;
\frac{(as;q)_\infty}{(az;q)_\infty}
\; \frac{(tq^{k})^{n-k}}{(atq^{k};q)_\infty}\; 
\notag \\  
&\qquad= \frac{(as;q)_\infty}{(az,at;q)_\infty}\sum_{k=0}^\infty 
\; \frac{\left(\frac{s}{z},at;q\right)_k\,z^k}{(as,q;q)_k}
\sum_{n=k}^\infty \frac{(r,f,g;q)_n\,u^n\, t^{n-k}}{(v,w,q;q)_n}\notag \\  
&\qquad=\frac{(as;q)_\infty}{(az,at;q)_\infty}\sum_{k=0}^\infty
\; \frac{\left(r,f,g,\frac{s}{z},at;q\right)_k\,(uz)^k}{(v,w,as,q;q)_k}
\sum_{n=0}^\infty \frac{(rq^k,fq^k,gq^k;q)_n(ut)^{n}}{(vq^k,wq^k,q;q)_n}.
\end{align}
Similarly, we have 
\begin{align}\label{proof3a}
&\mathbb{E}(r,f,g,v,w,u\theta_a)\left\{\frac{(az,at;q)_\infty}
{(as;q)_\infty }\right\}
\notag\\
&\quad =\sum_{n=0}^\infty \frac{(-1)^n\;q^{n\choose2}
(r,f,g;q)_nu^n}{(v,w,q;q)_n}
\theta_a^n\left\{\frac
{(az,at;q)_\infty}{(as;q)_\infty}\right\}\nonumber\\
&\quad=\sum_{n=0}^\infty \frac{(-1)^n\;q^{n\choose2}(r,f,g;q)_nu^n}{(v,w,q;q)_n}\; 
\sum_{k=0}^n \qbinomial{n}{k}{q}\theta_a^k\left\{\frac{(az;q)_\infty}
{(as;q)_\infty}\right\} \theta_a^{n-k}
\left\{(atq^{-k};q)_\infty \right\}\notag \\
&\quad=\sum_{n=0}^\infty \frac{(-1)^n\; q^{n\choose2}(r,f,g;q)_nu^n}
{(v,w,q;q)_n}\; 
\sum_{k=0}^n \qbinomial{n}{k}{q}
\frac{\left(-\frac{q}{a}\right)^k\; 
\left(\frac{z}{s};q\right)_k}{\left(\frac{q}{as};q\right)_k}\notag \\   
&\qquad \quad \cdot \frac{(az;q)_\infty}{(as;q)_\infty}
(atq^{-k};q)_\infty\; (-tq^{-k})^{n-k}\; 
\notag \\  
&\quad= \frac{(az,at;q)_\infty}{(as;q)_\infty}\sum_{k=0}^\infty 
\; \frac{(-1)^kq^{-{k\choose2}}\left(\frac{z}{s},
\frac{q}{at};q\right)_k\,t^k}
{\left(\frac{q}{as},q;q\right)_k}
\sum_{n=k}^\infty \frac{q^{{n\choose2}-k(n-k)}
(r,f,g;q)_n\,u^n\,t^{n-k}}
{(v,w,q;q)_n}\notag \\  
&\quad= \frac{(az,at;q)_\infty}{(as;q)_\infty}\sum_{k=0}^\infty 
\; \frac{\left(r,f,g,\frac{z}{s},\frac{q}{at};q\right)_k\,(-ut)^k}
{\left(v,w,\frac{q}{as},q;q\right)_k}
\sum_{n=0}^\infty \frac{q^{n\choose2}(rq^k,fq^k,gq^k;q)_n\,(ut)^{n}}
{(vq^k,wq^k,q;q)_n},
\end{align}
which evidently completes the proof of Theorem \ref{dsdpqro1}.
\end {proof}

We remark that, when $g=w=0$, Theorem \ref{dsdpqro1} reduces to 
the concluding result of Li and Tan \cite{LiTan2016}.

\begin{cor} \label{sdpqro1}
It is asserted that  
\begin{align}
\mathbb{T}(r,f,g,v,w,uD_s)\left\{ \frac{1 }{(xs;q)_\infty}\right\}
=\frac{1 }{(xs;q)_\infty} {}_3\Phi_2\left[\begin{matrix}
\begin{array}{rrr}
r,f,g;\\
\\
v,w;
\end{array}
\end{matrix} q;xu\right] 
\end{align}
and
\begin{align}
\mathbb{E}(r,f,g,v,w,-u\theta_s)\left\{(xs;q)_\infty\right\}
=(xs;q)_\infty\,{}_3\Phi_3\left[\begin{matrix}
\begin{array}{rrr}
r,f,g;\\
\\
v,w,0;
\end{array}
\end{matrix} q;xu\right],
\end{align}
provided that $\max\left\{|xs|,|xu|\right\}<1$.
\end{cor}

The goal in this paper is to give potentially 
useful generalizations of a number
$q$-series and $q$-integral identities 
such as the $q$-binomial 
theorem or the $q$-Gauss sum, 
the $q$-Chu-Vandermonde summation formula 
and the Andrews-Askey integral. 

Our paper is organized as follows. 
In Section \ref{generalize0}, we give two formal generalizations 
of the $q$-binomial theorem or the $q$-Gauss sum 
by applying the $q$-difference equations. 
In Section \ref{generalize1}, we derive a set of two extensions 
$q$-Chu-Vandermonde summation formulas by making use of the 
$q$-difference equations. Next, in Section \ref{generalize2}, 
we derive two new generalizations 
of the Andrews-Askey integral by means of the $q$-difference equations. 
Finally, in our last section (Section \ref{conclusion}), we present
a number of concluding remarks and observations concerning the various
results which we have considered in this investigation.
 
\section{\bf A Set of Formal Generalizations of the $q$-Binomial Theorem}
\label{generalize0}

We begin this section by recalling the following $q$-binomial theorem 
(see, for example, \cite{GasparRahman}, \cite{Slater} 
and \cite{SrivastavaKarlsson}):
\begin{equation}
\sum_{n=0}^\infty\frac{(a;q)_n\,x^n}{(q;q)_n}
=\frac{(ax;q)_\infty}{(x;q)_\infty}\qquad (|x|<1).\label{qbino}
\end{equation}
In Theorem \ref{thm_10} below, we give two 
generalizations of the $q$-binomial theorem 
\eqref{qbino} by applying the $q$-difference equations. 

\begin{thm}
\label{thm_10} 
Each of the following assertions holds true$:$
\begin{align}
&\sum_{n=0}^\infty\frac{(a;q)_n \;a^{-n}}{(q;q)_n}\;\sum_{k=0}^\infty\;
\frac{(q^{-n},ax;q)_k\,q^k}{(q;q)_k}\;\sum_{j,i \geqq 0}\;
\frac{(r,f,g;q)_{j+i}\;u^{j+i}}
{(q;q)_i \;(v,w;q)_{j+i}}\;\frac{\left(\frac{c}{b},axq^k;q\right)_{j}}
{(cx,q;q)_{j}}\;(aq^k)^i\;b^j\notag \\
&\qquad =\frac{(ax;q)_\infty }{(x;q)_\infty}\; \sum_{j,i\geqq 0}
\frac{(r,f,g;q)_{j+i}\;u^{j+i}}{(q;q)_i\; (v,w;q)_{j+i}}\;
\frac{\left(\frac{c}{b},x;q\right)_{j}}{(cx,q;q)_{j}}\; b^j \qquad 
(|x|<1)
\label{vlem10}
\end{align}
and
\begin{align}
&\sum_{n=0}^\infty\;\frac{(a;q)_n \;a^{-n}}{(q;q)_n}\;\sum_{k=0}^\infty\;
\frac{(q^{-n},ax;q)_k\,q^k}{(q;q)_k}\;\sum_{i,j\geqq 0} \; 
\frac{(-1)^{j+i}\;q^{\left({}^i_2\right)}(r,f,g;q)_{j+i}}
{(q;q)_i\;(v,w;q)_{j+i}}\notag \\
&\qquad \qquad \qquad \cdot \frac{\left(\frac{bq^{1-k}}{a},
\frac{q}{cx};q\right)_j}{\left(\frac{q^{1-k}}{ax},q;q\right)_j}
\,(uc)^{j+i} \notag \\
& \qquad= \frac{(ax;q)_\infty}{(x;q)_\infty}\; \sum_{i,j=0}^\infty 
\; \frac{(-1)^{j+i}\;q^{\left({}^i_2\right)}\;(r,f,g;q)_{j+i}}
{(q;q)_i(v,w;q)_{j+i}}\; \frac{\left(\frac{q}{bx};q\right)_j}
{\left(\frac{q}{x},q;q\right)_j} \,(bu)^{j+i}, \label{thlem10}
\end{align}
provided that both sides of $\eqref{vlem10}$ and $\eqref{thlem10}$ exist.
\end{thm}

\begin{rem}
{\rm For $u=0$ and by using the fact that  
$$\sum_{k=0}^\infty\frac{(q^{-n},ax;q)_k}{(q;q)_{k}}\;q^k
={}_2\Phi_1\left[\begin{matrix}
\begin{array}{rrr}
q^{-n}, ax;\\
\\
0;
\end{array}
\end{matrix} q; q\right]=(ax)^n,
$$
the assertions $\eqref{vlem10}$ and $\eqref{thlem10}$ reduce  
to $\eqref{qbino}$.}
\end{rem}

In our proof of Theorem \ref{thm_10}, we shall need Theorem \ref{thm_121} 
and Corollary \ref{cor_11} below.
 
\begin{thm}
\label{thm_121}
Each of the following assertions holds true$:$
\begin{align}
&\mathbb{T}(r,f,g,v,w,uD_x)\left\{\frac{p_n\left(x,\frac{y}{a};q\right)_n\;
(cx;q)_\infty}{(ax,bx;q)_\infty}\right\} \notag \\
&\qquad =\frac{(y;q)_n}{a^n}\; \frac{(cx;q)_\infty}{(ax,bx;q)_\infty}
\sum_{k=0}^\infty\frac{(q^{-n},ax;q)_k\,q^k}{(y,q;q)_k}\notag \\
&\qquad \qquad \cdot \sum_{j,i\geqq 0}\frac{(r,f,g;q)_{j+i}\;u^{j+i}}
{(q;q)_i (v,w;q)_{j+i}}\; 
\frac{\left(\frac{c}{b},axq^k;q\right)_{j}}{(cx,q;q)_{j}}\;
(aq^k)^i\;b^j\label{thm10}
\end{align}
and
 
\begin{align}
&\mathbb{E}(r,f,g,v,w,u\theta_x)\left\{\frac{p_n\left(x,\frac{y}{a}\right)
(bx,cx;q)_\infty}{(ax;q)_\infty}\right\}\notag \\
&\qquad =\frac{(y;q)_n }{a^ n }\frac{(bx,cx;q)_\infty}{(ax;q)_\infty}\;\sum_{i,j\geqq 0}  
\frac{(-1)^{j+i}\;q^{({}^i_2)}\;(r,f,g;q)_{j+i}}{(q;q)_i\;(v,w;q)_{j+i}}\notag \\
&\qquad \qquad \cdot \frac{\left(\frac{bq^{1-k}}{a},\frac{q}{cx};q\right)_j}
{\left(\frac{q^{1-k}}{ax},q;q\right)_j}
\,(uc)^{j+i},\label{2f2.16}
\end{align}
provided that $\max\left\{|ax|,|bx|,|cx|\right\}<1$.
\end{thm}

\begin{cor}
\label{cor_11} 
Each of the following assertions holds true$:$
\begin{align}
&\mathbb{T}(r,f,g,v,w,uD_x)\left\{\frac{x^n(cx;q)_\infty}
{(ax,bx;q)_\infty}\right\}\notag \\
&\qquad=\frac{1}{a^n}\;\frac{(cx;q)_\infty}{(ax,bx;q)_\infty}\;\sum_{k=0}^\infty
\; \frac{(q^{-n},ax;q)_k\,q^k}{(q;q)_k}\notag \\
&\qquad \qquad \cdot \sum_{j,i\geqq 0}\;\frac{(r,f,g;q)_{j+i}\;u^{j+i}}
{(q;q)_i \;(v,w;q)_{j+i}}\; 
\frac{\left(\frac{c}{b},axq^k;\right)_{j}}{(cx,q;q)_{j}}\;(aq^k)^i\;b^j 
\label{lem10}
\end{align}
and  
 
\begin{align}
&\mathbb{E}(r,f,g,v,w,u\theta_x)\left\{\frac{x^n(cx,bx;q)_\infty}
{(ax;q)_\infty}\right\}\notag \\ 
&\qquad = \frac{1 }{a^ n }\frac{(bx,cx;q)_\infty}{(ax;q)_\infty}\sum_{k=0}^\infty\;
\frac{(q^{-n},ax;q)_k\,q^k}{(q;q)_k}\sum_{i,j\geqq 0}  
\frac{(-1)^{j+i}q^{({}^i_2)}\;(r,f,g;q)_{j+i} }{(q;q)_i\;(v,w;q)_{j+i}}
\notag \\
&\qquad \qquad \cdot \frac{\left(\frac{bq^{1-k}}{a},\frac{q}{cx};q\right)_j}
{\left(\frac{q^{1-k}}{ax},q;q\right)_j}
\;(cu)^{j+i}, \label{ff2.16}
\end{align}
provided that $\max\left\{|ax|,|bx|,|cx|,|cu|\right\}<1$.
\end{cor}

\begin{rem}
{\rm For $y=0,$ the assertions $\eqref{thm10}$ and $\eqref{2f2.16}$ 
reduce  to $\eqref{lem10}$ and $\eqref{ff2.16},$ respectively.}
\end{rem}

\begin{proof}[Proof of Theorem $\ref{thm_121}$]
Upon first setting  $x\to ax$ in (\ref{ttf}) and then multiplying 
both sides of the resulting equation by 
$\frac{(cx;q)_\infty}{(bx;q)_\infty},$ we get
\begin{equation}
\label{9z}
\sum_{n=0}^\infty\;\frac{(q^{-n};q)_k\,q^k}{(y,q;q)_k}\;
\frac{(cx;q)_\infty}{(axq^k,bx;q)_\infty}
=\frac{(ax)^n\left(\frac{y}{ax};q\right)_n\;
(cx;q)_\infty}{(y;q)_n(ax,bx;q)_\infty}.
\end{equation}
Now, by applying the operator $\mathbb{T}(r,f,g,v,w,uD_x)$ 
to both sides of (\ref{9z}), it is easy to see that 
\begin{align}
&\sum_{k=0}^\infty\frac{(q^{-n};q)_k\,q^k}{(y,q;q)_k}\;
\mathbb{T}(r,f,g,v,w,uD_x)
\left\{\frac{(cx;q)_\infty}{(axq^k,bx;q)_\infty}\right\}
\notag \\
&\qquad=\frac{a^ n}{(y;q)_n}\;
\mathbb{T}(r,f,g,v,w,uD_x)
\left\{\frac{x^n\left(\frac{y}{ax};q\right)_n(cx;q)_\infty}
{(ax,bx;q)_\infty}\right\}\notag \\
&\qquad= \frac{a^n}{(y;q)_n}\;\mathbb{T}(r,f,g,v,w,uD_x)
\left\{\frac{p_n\left(x,\frac{y}{a}\right)(cx;q)_\infty}
{(ax,bx;q)_\infty}\right\}.\label{2.16}
\end{align}

The proof of the first assertion (\ref{thm10}) of Theorem \ref{thm_121} 
is completed by using the relation (\ref{gLEM})
in the left-hand side of (\ref{2.16}). \\
The proof of the second assertion (\ref{2f2.16}) of Theorem \ref{thm_121}  
is much akin to that of the first assertion (\ref{thm10}). The details 
involved are, therefore, being omitted here. 
\end{proof}

\begin{proof}[Proof of Theorem $\ref{thm_10}$]
Multiplying both sides of (\ref{qbino}) by 
$\frac{(cx;q)_\infty}{(bx;q)_\infty}$, we find that
\begin{equation}
\sum_{n=0}^\infty\;\frac{(a;q)_n}{(q;q)_n}\; 
\frac{x^n\,(cx;q)_\infty}{(ax,bx;q)_\infty} 
=\frac{(cx;q)_\infty}{(bx,x;q)_\infty}.\label{2qbino}
\end{equation}
Eq. (\ref{vlem10}) can be written equivalently as follows:
\begin{align}
&\sum_{n=0}^\infty\;\frac{(a;q)_n}{(q;q)_n}\cdot
\frac{a^{-n}(cx;q)_\infty}{(bx,ax;q)_\infty} 
\;\sum_{k=0}^\infty\;\frac{(q^{-n},ax;q)_k\,q^k}{(q;q)_k}\notag \\
&\qquad \qquad \qquad \cdot \sum_{j,i\geqq 0}\;\frac{(r,f,g;q)_{j+i}\;
u^{j+i}}{(q;q)_i (v,w;q)_{j+i}} 
\frac{\left(\frac{c}{b},axq^k;q\right)_{j}}{(cx,q;q)_{j}}\;(aq^k)^i\;b^j
\notag \\
& \qquad=\frac{(cx;q)_\infty}{(bx,x;q)_\infty}\;
\sum_{j,i=0}\frac{(r,f,g;q)_{j+i}\;u^{j+i}}{(q;q)_i (v,w;q)_{j+i}} 
\cdot \frac{\left(\frac{c}{b},x;q\right)_{j}}{(cx,q;q)_{j}}\; b^j . 
\label{rvlem10}
\end{align}

If we use $F(r,f,g,v,w,a,u)$ to denote the right-hand side 
of (\ref{rvlem10}),
it is easy to verify that $F(r,f,g,v,w,a,u)$ satisfies  
(\ref{thm2_1}).  By applying (\ref{thm2_1.2}), we thus find that
\begin{align}
\label{efss}
F(r,f,g,v,w,x,u)&=\mathbb{T}(r,f,g,v,w,uD_x)\Big\{F(r,f,g,v,w,x,0)\Big\} 
\notag \\
&=\mathbb{T}(r,f,g,v,w,uD_x)\left\{\frac{(cx;q)_\infty}{(bx,x;q)_\infty}\right\}
\notag \\
&=\mathbb{T}(r,f,g,v,w,uD_x)\left\{\sum_{n=0}^\infty\;\frac{(a;q)_n}
{(q;q)_n}\; \frac{x^n\,(cx;q)_\infty}{(ax,bx;q)_\infty} \right\}
\notag \\
&=\sum_{n=0}^\infty\;\frac{(a;q)_n}{(q;q)_n}\;  
\mathbb{T}(r,f,g,v,w,uD_x)\left\{\frac{x^n\,(cx;q)_\infty}{(ax,bx;q)_\infty}  
\right\}.
\end{align}

The proof of the first assertion (\ref{vlem10}) of Theorem \ref{thm_10}  
can now be completed by making use of the relation (\ref{lem10}).

The proof of the second assertion (\ref{thlem10}) of Theorem \ref{thm_10}  
is much akin to that of the first assertion (\ref{vlem10}). The details 
involved are, therefore, being omitted here. 
\end{proof}

\section{\bf Two Generalizations of the $q$-Chu-Vandermonde Summation Formula}
\label{generalize1}

The $q$-Chu-Vandermonde summation formula is recalled here as follows 
(see, for example, \cite{GEAndrews1986} and \cite{GasparRahman}): 
\begin{equation} \label{ttf}
{}_2\Phi_1\left[\begin{matrix}
\begin{array}{rrr}
q^{-n},x;\\
\\
y;
\end{array}
\end{matrix} q; q\right]=\frac{\left(\frac{y}{x};q\right)_n}
{(y;q)_n}\;x^n
\qquad \big(n\in \mathbb{N}_0:=\mathbb{N}\cup\{0\}\big).
\end{equation}
 
In this section, we give two generalizations of the
$q$-Chu-Vandermonde summation formula \eqref{ttf} 
by applying $q$-difference equations. 

\begin{thm}
\label{thm_11}
The following assertion holds true for $y \neq 0$$:$ 
\begin{align}
&\sum_{k=0}^{n} \frac{(q^{-n},x;q)_k\,q^k}{(q,y;q)_k}
\;{}_3\Phi_2\left[\begin{matrix}
\begin{array}{rrr}
r,f,g;\\
\\
v,w;
\end{array}
\end{matrix} q; uq^k\right] \notag \\
&\qquad=\frac{x^n\left(\frac{y}{x};q\right)_n}{(y;q)_n}\; \sum_{k,j\geqq 0}
\frac{(r,f,g;q)_{k+j}}{(q;q)_{j}\;(v,w;q)_{k+j}}\;
\frac{\left(\frac{q^{1-n}}{y},\frac{qx}{y};q\right)_{k}}
{\left(\frac{xq^{1-n}}{y},q;q\right)_k}\; u^{k+j}\left(\frac{q}{y}\right)^j.
\label{qCh}
\end{align}
\end{thm}

We next derive another generalization of the $q$-Chu-Vandermonde 
summation formula \eqref{ttf} as follows.

\begin{thm}
\label{thm_1f1} 
For $m\in\mathbb{N}_0$ and $y\neq 0$$,$ 
it is asserted that
\begin{align}
{}_2\Phi_1\left[\begin{matrix}
\begin{array}{rrr}
q^{-n},x;\\
\\
y;
\end{array}
\end{matrix} q; q^{1+m}\right] 
=\frac{x^n\left(\frac{y}{x};q\right)_n}{(y;q)_n}   
\;\sum_{j= 0}^m\begin{bmatrix}
m \\
j \\
\end{bmatrix}_q \;\frac{\left(\frac{q^{1-n}}{y},\frac{qx}{y};q\right)_{m-j}}
{\left(\frac{xq^{1-n}}{y};q\right)_{m-j}} \; \left(\frac{q}{y}\right)^j. 
\label{AqCh}
\end{align}
\end{thm}

\begin{rem}
{\rm For $u=0$ or $m=0,$ the assertion $\eqref{qCh}$ or $\eqref{AqCh}$ reduces to
the $q$-Chu-Vandermonde summation formula $\eqref{ttf}$.  
Furthermore, if we first set $i+j=m$ and then extract the coefficients 
of $\displaystyle \frac{(r,f,g;q)_m}{(v,w;q)_m}u^m$ 
from the two members of the assertion $\eqref{qCh}$ of Theorem $\ref{thm_11},$  
we obtain the transformation formula (\ref{AqCh}), which leads us to 
the $q$-Chu-Vandermonde summation formula $\eqref{ttf}$ when $m=0$. 
Also, upon putting $n=0,$ the assertion $\eqref{AqCh}$ 
reduces to the following identity$:$}
\begin{align}
\sum_{j=0}^m \begin{bmatrix}
m \\
j \\
\end{bmatrix}_q\;\left(\frac{q}{y};q\right)_{m-j}\;
\left(\frac{q}{y}\right)^j= 1\qquad (y\neq 0).
\end{align}
\end{rem}

\begin{proof}[Proof of Theorem $\ref{thm_11}$] 
We first write (\ref{ttf}) in the following form:
\begin{equation}
\sum_{k=0}^{n}\; \frac{(q^{-n};q)_k\,q^k}{(y,q;q)_k}
\;\frac{1}{(xq^k;q)_\infty}
=\frac{(-1)^{n}\;y^n\; q^{\left({}^{n}_{2}\right)}}
{(y;q)_n}\frac{\left(\frac{xq^{1-n}}{y};q\right)_\infty}
{\left(x,\frac{qx}{y};q\right)_\infty}.
\end{equation}
Eq. (\ref{qCh}) can be written equivalently as follows:
\begin{align}
&\sum_{k=0}^\infty\;\frac{(q^{-n};q)_k\,q^k}{(q,c;q)_k}\cdot
\frac{1}{(xq^k;q)_\infty}\; {}_3\Phi_2\left[\begin{matrix}
\begin{array}{rrr}
r,f,g;\\
\\
v,w;
\end{array}
\end{matrix} q; uq^k\right] \notag \\
&\qquad=\frac{(-1)^{n}\;y^n \;q^{\left({}^{n}_{2}\right)}}
{(y;q)_n}\;\frac{\left(\frac{xq^{1-n}}{y};q\right)_\infty}
{\left(x,\frac{qx}{y};q\right)_\infty}\cdot\sum_{i,j\geqq 0}\; 
\frac{(r,f,g;q)_{j+i}\;u^{j+i}}{(q;q)_i\;(v,w;q)_{j+i}}\;
\frac{\left(\frac{q^{1-n}}{y},\frac{qx}{y};q\right)_{j}}
{\left(\frac{xq^{1-n}}{y},q;q\right)_{j}}\;\left(\frac{q}{y}\right)^i.
\label{qdqdqq}
\end{align}

If we use $G(r,f,g,v,w,x,u)$ to denote the right-hand side of (\ref{qdqdqq}),
it is easy to observe that $G(r,f,g,v,w,x,u)$ satisfies (\ref{thm2_1}). By  
using (\ref{thm2_1.2}), we obtain
\begin{align}
\label{efss1}
G(r,f,g,v,w,x,u)&=\mathbb{T}(r,f,g,v,w,uD_x)\Big\{G(r,f,g,v,w,x,0) \Big\} 
\notag \\
&=\mathbb{T}(r,f,g,v,w,uD_x)\left\{\frac{(-1)^{n}\;y^n\; q^{\left({}^{n}_{2}\right)}}{(y;q)_n}
\;\frac{\left(\frac{xq^{1-n}}{y};q\right)_\infty}{\left(x,\frac{qx}{y};q\right)_\infty}\right\}
\notag \\
&=\mathbb{T}(r,f,g,v,w,uD_x)\left\{\sum_{k=0}^\infty\frac{(q^{-n};q)_k\,q^k}
{(y,q;q)_k}\; \frac{1}{(xq^k;q)_\infty} \right\}
\notag \\
&=\sum_{k=0}^{n}\frac{(q^{-n};q)_k\,q^k}{(y,q;q)_k}  \;  
\mathbb{T}(r,f,g,v,w,uD_x)\left\{\frac{1}{(xq^k;q)_\infty}   
\right\}.
\end{align}

Finally, by using the fact that
\begin{equation}
\mathbb{T}(r,f,g,v,w,uD_x)\left\{\frac{1}{(xq^k;q)_\infty}\right\} 
=\frac{1}{(xq^k;q)_\infty} \;
{}_3\Phi_2\left[\begin{matrix}
\begin{array}{rrr}
r,f,g;\\
\\
v,w;
\end{array}
\end{matrix} 
q; uq^k\right], \label{fgLEM}   
\end{equation}
and after some simplification involving $\frac{1}{(x;q)_\infty}$, 
we get the left-hand side of (\ref{qCh}).
\end{proof}

\section{\bf New Generalizations of the Andrews-Askey Integral}
\label{generalize2}

The following famous formula is known as the 
Andrews-Askey integral (see, for details, \cite{GEA-RA1981}). 
It was derived from Ramanujan's celebrated
${}_1\Psi_1$-summation formula.
 
\begin{pro}{\rm (see \cite[Eq. (2.1)]{GEA-RA1981})}. 
For $\max\left\{|ac|,|ad|,|bc|,|bd|\right\}<1,$
it is asserted that
\begin{align}
\label{eqd}
\int_c^d\frac{\left(\frac{qt}{c},\frac{qt}{d};q\right)_\infty}
{(at,bt;q)_\infty}\; {\rm d}_qt
=\frac{d(1-q)\left(q,\frac{dq}{c},\frac{c}{d},abcd;q\right)_\infty}
{(ac,ad,bc,bd;q)_\infty}.
\end{align}
\end{pro}

The Andrews-Askey integral \eqref{eqd} is indeed an important 
formula in the theory of $q$-series (see \cite{Liu97}). 
  
Recently, Cao \cite{JianCao2013} gave the following two generalizations of 
the Andrews-Askey integral \eqref{eqd} by the method based upon 
$q$-difference equations. 

\begin{pro}{\rm (see \cite[Theorems 14 and 15]{JianCao2013})}
\label{aznzzT} 
For $N\in\mathbb{N}$ and $r=q^{-N},$ suppose that 
$$\max\left\{|ac|,|ad|,|bc|,|bd|,\left|\frac{qwr}{v}\right|,
\left|\frac{q}{v}\right|\right\}<1.$$ 
Then
\begin{align}
\label{hdqd}
&\int_c^d\; \frac{\left(\frac{qt}{c},\frac{qt}{d};q\right)_\infty}
{(at,bt;q)_\infty}  
\, {}_4\Phi_2\left[\begin{matrix}
\begin{array}{rrr}
r,w,\frac{c}{t},abcd;\\
\\
ac,\frac{qwr}{v};
\end{array}
\end{matrix} q; \frac{qt}{vbcd}\right]\; {\rm d}_qt\notag \\
&\qquad \quad =\frac{d(1-q)\left(q,\frac{dq}{c},\frac{c}{d},abcd,
\frac{qw}{v},\frac{qr}{v};q\right)_\infty}
{\left(ac,ad,bc,bd,\frac{qwr}{v},\frac{q}{v};q\right)_\infty} 
\; {}_2\Phi_1\left[\begin{matrix}
\begin{array}{rrr}
w,r;\\
\\
v;
\end{array}
\end{matrix} q; \frac{q}{bc}\right].
\end{align}
Furthermore$,$ for $N\in\mathbb{N}$ and $r=q^{-N},$ suppose that 
$$\max\left\{|ac|,|ad|,|bc|,|bd|,\left|\frac{v}{w}\right|,
\left|\frac{v}{r}\right|\right\}<1.$$ 
Then
\begin{align}
\label{1hdqd}
&\int_c^d\;\frac{\left(\frac{qt}{c},\frac{qt}{d};q\right)_\infty}
{(at,bt;q)_\infty}  
\, {}_4\Phi_2\left[\begin{matrix}
\begin{array}{rrr}
r,w,\frac{c}{t},\frac{q}{ad};\\
\\
\frac{q}{at},\frac{qrw}{v};
\end{array}
\end{matrix} q;  q\right]\; d_qt \notag \\
&\qquad \quad=\frac{d(1-q)\left(q,\frac{dq}{c},
\frac{c}{d},abcd,\frac{v}{wr},v; q\right)_\infty}
{\left(ac, ad, bc, bd,\frac{v}{w},\frac{v}{r};q\right)_\infty} 
\; {}_2\Phi_1\left[\begin{matrix}
\begin{array}{rrr}
w,r;\\
\\
v;
\end{array}
\end{matrix} q; \frac{vbc}{wr}\right].
\end{align}
\end{pro}

In this section, we give the following two generalizations of the 
Andrews-Askey integral \eqref{eqd} by using the method 
of $q$-difference equations.
    
\begin{thm}
\label{fgazzzT} 
For $M\in\mathbb{N}$ and $r=q^{-M},$ suppose that 
$$\max\left\{|ac|,|ad|,|bc|,|bd|,\left|\frac{q}{bc}\right|\right\}<1.$$ 
Then
\begin{align}
\label{gqdqd}
&\int_c^d\;\frac{\left(\frac{qt}{c},\frac{qt}{d};q\right)_\infty}
{(at,bt;q)_\infty} 
\;\sum_{k=0}^\infty
\; \frac{\left(r,f,g,\frac{c}{t},abcd;q\right)_{k}\;
\left(\frac{qt}{bcd}\right)^k}
{(v,w,ac,q;q)_{k}}\; {}_3\Phi_2\left[\begin{matrix}
\begin{array}{rrr}
rq^k,fq^k,gq^k;\\
\\
vq^k,wq^k;
\end{array}
\end{matrix} q; q\right] \; {\rm d}_qt \notag \\
&\qquad \quad =\frac{d(1-q)\left(q,\frac{dq}{c},
\frac{c}{d},abcd; q\right)_\infty}
{(ac, ad, bc, bd;q)_\infty} \; 
{}_3\Phi_2\left[\begin{matrix}
\begin{array}{rrr}
r,f,g;\\
\\
v,w;
\end{array}
\end{matrix} q; \frac{q}{bc}\right].
\end{align}
\end{thm}

\begin{thm}
\label{vT} 
For $M\in\mathbb{N}$ and $r=q^{-M},$ suppose that 
$\max\left\{|ac|,|ad|,|bc|,|bd|\right\}<1$.
Then
\begin{align}
\label{vdqd}
&\int_c^d\;\frac{\left(\frac{qt}{c},\frac{qt}{d};q\right)_\infty} 
{(at,bt;q)_\infty}   
\;\sum_{k=0}^\infty 
\; \frac{\left(r,f,g,\frac{c}{t},\frac{q}{ad};q\right)_k\,
\left(\frac{vw}{rfg}\right)^k}
{\left(v,w,\frac{q}{at},q;q\right)_k}
\; {}_3\Phi_3\left[\begin{matrix}
\begin{array}{rrr}
rq^k,fq^k,gq^k;\\
\\
vq^k,wq^k,0;
\end{array}
\end{matrix} q;-\frac{vw}{rfg}\right]\; {\rm d}_qt \notag \\
& \qquad \quad =\frac{d(1-q)\left(q,\frac{dq}{c},
\frac{c}{d},abcd; q\right)_\infty}
{(ac, ad, bc, bd;q)_\infty}\;{}_3\Phi_3\left[\begin{matrix}
\begin{array}{rrr}
r,f,g;\\
\\
v,w,0;
\end{array}
\end{matrix} q; \frac{vwbc}{rfg}\right].
\end{align}
\end{thm}

\begin{rem}
{\rm For $r=1,$ both $\eqref{gqdqd}$ and $\eqref{vdqd}$ 
reduce  to $\eqref{eqd}$.  
Moreover$,$ for $r=q^{-N},$ $g=w=0$ and $u=\frac{q}{bcd},$ 
the assertion (\ref{gqdqd}) of 
Theorem \ref{fgazzzT} reduces to $\eqref{hdqd}$. 
For $r=q^{-N},$ $g=w=0$ and $u=\frac{v}{rfbcd},$ 
the assertion $\eqref{gqdqd}$ of  
Theorem $\ref{vT}$ reduces to $\eqref{vdqd}$.} 
\end{rem}
 
\begin{proof}[Proof of Theorems $\ref{fgazzzT}$ and $\ref{vT}$]
Eq. (\ref{gqdqd}) can be written equivalently as follows:
\begin{align}
\label{qdqdqq1}
&\int_c^d\;\frac{\left(\frac{qt}{c},\frac{qt}{d};q\right)_\infty} 
{(bt;q)_\infty} \cdot \frac{(ac;q)_\infty}{(at, abcd;q)_\infty}\; 
\sum_{k=0}^\infty\;\frac{\left(r,f,g,\frac{c}{t},abcd;q\right)_{k}
\left(\frac{qt}{bcd}\right)^k}{(v,w,ac,q;q)_{k}}\notag \\
&\qquad \qquad \qquad \cdot{}_3\Phi_2\left[\begin{matrix}
\begin{array}{rrr}
rq^k,fq^k,gq^k;\\
\\
vq^k,wq^k;
\end{array}
\end{matrix} q; q\right] \; {\rm d}_qt\notag \\
&\qquad \quad = \frac{d(1-q)\left(q,\frac{dq}{c},
\frac{c}{d};q\right)_\infty}
{(bc, bd;q)_\infty}\cdot\frac{1}{(ad;q)_\infty}\;
{}_3\Phi_2\left[\begin{matrix}
\begin{array}{rrr}
r,f,g;\\
\\
v,w;
\end{array}
\end{matrix} q; \frac{q}{bc}\right].
\end{align}

If we use $H(r,f,g,v,w,a,u)$ to denote the right-hand side of (\ref{qdqdqq1}),
it is easy to see that $H(r,f,g,v,w,a,u)$ satisfies (\ref{thm2_1}) 
with $u=\frac{q}{bcd}$. By making use of 
(\ref{thm2_1.2}), we thus find that
\begin{align*}
H(r,f,g,v,w,a,u)&=\mathbb{T}(r,f,g,v,w,uD_a)\Big\{H(r,f,g,v,w,a,0)\Big\}
\notag \\
&=\mathbb{T}\left(r,f,g,v,w,\frac{q}{bcd}\;D_a\right)
\Big\{H(1,f,g,v,w,a,u) \Big\} \notag \\
&=\mathbb{T}\left(r,f,g,v,w,\frac{q}{bcd}\;D_a\right)
\left\{\frac{d(1-q)\left(q,\frac{dq}{c},\frac{c}{d}; q\right)_\infty}
{(bc, bd;q)_\infty}\cdot\frac{1}{(ad;q)_\infty} \right\}
\notag \\
&=\mathbb{T}\left(r,f,g,v,w,\frac{q}{bcd}\;D_a\right)
\left\{\int_c^d\;\frac{\left(\frac{qt}{c},\frac{qt}{d};q\right)_\infty} 
{(bt;q)_\infty} \cdot \frac{(ac;q)_\infty}{(at, abcd;q)_\infty}\; 
{\rm d}_qt \right\}
\notag \\
&=\int_c^d\; \frac{\left(\frac{qt}{c},\frac{qt}{d};q\right)_\infty} 
{(bt;q)_\infty} \cdot  
\mathbb{T}\left(r,f,g,v,w,\frac{q}{bcd}\;D_a\right)
\left\{\frac{(ac;q)_\infty}{(at,abcd;q)_\infty}  
\right\} {\rm d}qt.
\end{align*}
Now, by applying the fact that
\begin{align*}
&\mathbb{T}\left(r,f,g,v,w,\frac{q}{bcd}\;D_a\right)
\left\{\frac{(ac;q)_\infty}{(at,abcd;q)_\infty}\right\}
\notag \\
&\qquad \quad=\frac{(ac;q)_\infty}{(at,abcd;q)_\infty}\; 
\sum_{k=0}^\infty\;
\frac{\left(r,f,g,\frac{c}{t},abcd;q\right)_{k}\;
\left(\frac{qt}{bcd}\right)^k}
{(v,w,ac,q;q)_{k}}\; {}_3\Phi_2\left[\begin{matrix}
\begin{array}{rrr}
rq^k,fq^k,gq^k;\\
\\
vq^k,wq^k;
\end{array}
\end{matrix} q; q\right],
\end{align*}
we get the left-hand side of (\ref{gqdqd}).

The proof of the assertion (\ref{vdqd}) of Theorem \ref{vT}  
is much akin to that of the assertion (\ref{gqdqd})
of Theorem \ref{fgazzzT}. The details involved are, 
therefore, being omitted here.
 
The proofs of Theorems \ref{fgazzzT} and \ref{vT} 
are thus completed.
\end{proof}

\section{\bf Concluding Remarks and Observations}
\label{conclusion}

In our present investigation, we have introduced a set of 
two $q$-operators \break $\mathbb{T}(a,b,c,d,e,yD_x)$ and 
$\mathbb{E}(a,b,c,d,e,y\theta_x)$ with to applying them
to derive two potentially useful 
generalizations of the $q$-binomial theorem, two extensions of 
the $q$-Chu-Vandermonde summation formula and two new generalizations 
of the Andrews-Askey integral by means of the $q$-difference equations.
We have also briefly described relevant 
connections of various special cases 
and consequences of our main results 
with several known results.

It is believed that the $q$-series and $q$-integral identities, which we
have presented in this paper, as well as the various related recent works
cited here, will provide encouragement and motivation for further researches
on the topics that are dealt with and investigated in this paper.\\

\medskip 

\noindent
{\bf Conflicts of Interest:} The authors declare that they 
have no conflicts of interest.\\

\end{document}